\documentclass[12pt]{iopart}
\usepackage{amssymb}
\usepackage{graphicx}
\begin{document}

\setcounter{section}{0}
\newtheorem{theorem}{Theorem}
\newtheorem{lemma}{Lemma}
\newtheorem{remark}{Remark}
\newtheorem{proposition}{Proposition}
\newtheorem{corollary}{Corollary}

\eqnobysec

\title[Formation of singularities in solutions to ideal hydrodynamics of inelastic gases]
{Formation of singularities in solutions to ideal hydrodynamics of
freely cooling inelastic gases}

\author{ Olga Rozanova}

\address{Department of Differential Equations \& Mechanics and Mathematics Faculty,
Moscow State University, Moscow, 119992,
 Russia}
\ead{rozanova@mech.math.msu.su}

\begin{abstract}

We consider  solutions to the hyperbolic system of equations of
ideal granular hydrodynamics  with conserved mass, total energy and
finite momentum of inertia and prove that these solutions
generically lose the initial smoothness within a finite time in any
space dimension $n$  for the adiabatic index $\gamma \le
1+\frac{2}{n}.$ Further, in the one-dimensional case we introduce a
solution  depending only on  the spatial coordinate outside of a
ball containing the origin and  prove that this solution under
rather general assumptions on initial data cannot be global in time
too. Then we construct an exact  axially  symmetric solution with
separable time and space variables having a strong singularity in
the density component beginning from the initial moment of time,
whereas other components of solution are initially continuous.
\end{abstract}

\ams{35L60,\, 76N10, \,35L67}

\maketitle

\section{Introduction}
\label{1}

The motion of the dilute gas  where the characteristic hydrodynamic
length scale of the flow  is sufficiently large and the viscous and
heat conduction terms can be neglected  is governed by the systems
of equations of ideal granular hydrodynamics \cite{Brilliantov}.

This system is given in ${\mathbb R}\times{\mathbb R}^n,\, n\ge 1,$
and has the following form:
\begin{equation}\label{1.1}
\partial_t \rho+{\rm div}_x (\rho u)=0,
\end{equation}
\begin{equation}\label{1.2}
\partial_t(\rho u)+{\rm Div}_x (\rho u \otimes u)=-\nabla_x
p,
\end{equation}
\begin{equation}\label{1.3}
\partial_t T+(u,\nabla_x T) + (\gamma-1) T {\rm div}_x u=-\Lambda
\rho T^{3/2},
\end{equation}
where $\rho$ is the gas density, $ u=(u_1,...,u_n)$ is the velocity,
$T$ is the temperature, $p=\rho T$ is the pressure, and $\gamma$ is
the adiabatic index $(1<\gamma\le 1+\frac{2}{n}),$
$\Lambda=const>0.$ we denote ${\rm Div}_x$ and ${\rm div}_x$ the
divergency of tensor and vector, respectively, with respect to the
space variables. The only difference between equations
\eref{1.1}--\eref{1.3} and the standard ideal gas dynamic equations
(where the elastic colliding of particles is supposed) is the
presence  of the inelastic energy loss term $-\Lambda \rho T^{3/2}$
in \eref{1.3}.

The granular gases are now popular subject of experimental,
numerical and theoretical investigation (e.g. \cite{Brilliantov},
\cite{Ludvig}, \cite{Meerson2} and references therein). In contrast
to ordinary molecular gases, granular gases "cool" spontaneously
because of the inelastic collisions between the particles. The
inelasticity of the collisions generally causes the granular gas to
form dense clusters. The formation of complex structure of clusters
has been investigated by means of molecular dynamics simulations and
hydrodynamic simulations.

 The Navier-Stokes granular hydrodynamics is the natural language for
 a theoretical description of granular macroscopic flows.  A characteristic feature
of  time-dependent solutions of the continuum equations is a
formation of finite-time singularities: the density blowup signals
the formation of close-packed clusters.

System \eref{1.1} -- \eref{1.3} can be written in a  hyperbolic
symmetric form in variables $\,\rho, \,u, K=p\rho^{-\gamma}\,$ and
therefore the Cauchy problem
$$
\rho, u, p\rho^{-\gamma}\big|_{t=0}\in H^m({\mathbb R}^n), \,m\ge
1+\left[ \frac{n}{2}\right] .
$$
has a solution as smooth as initial data at least for small $t>0$
\cite{Kato}.

We will call the solution to \eref{1.1} -- \eref{1.3} {\it
classical} if $\rho>0, p>0$ and the components of solution belongs
to $C^1([0,T), H^m({\mathbb R}^n)),\,T\le\infty.$

System \eref{1.1} -- \eref{1.3} has no constant solution except the
trivial one $(p\equiv 0),$ therefore the solution with the
components highly decreasing as $|x|\to\infty$ can be considered as
a natural perturbation if this steady state in the case of the mass
conservation.

Let us note that there exists a solution (the homogeneous cooling
state) with constant $\rho, u$ and $p\ne 0$. In this case the
temperature $T =T(t)=(\frac{\Lambda \rho_0 t}{2}+T(0)^{-1})^{-2},$
where $T(0)$ is the initial value of temperature (the Haff's law).
Another trivial solution is $u=p=T\equiv 0, \, \rho(t,x)=\rho_0(x)$.

 We introduce the following
integrals: the total mass
$$M(t)=\int\limits_{{\mathbb R}^n}\rho \, d x,$$  the momentum
$$P(t)=\int\limits_{{\mathbb R}^n}\rho u\, d x,$$ and  the total energy
$$\mathcal E(t)=\int\limits_{{\mathbb R}^n}\left(\frac{1}{2}\rho |u|^2+\rho
T\right)\, d x \,=E_{k}(t)+E_{i}(t).$$ Here $E_k(t)$ and $E_i(t)$
are the kinetic and internal components of energy, respectively. Let
us introduce also the  functionals
$$G(t)=\frac{1}{2} \int\limits_{{\mathbb
R}^n}\rho(t,x)|{ x}|^2\,dx,\qquad F(t)=\int\limits_{{\mathbb R}^n}
({u},{ x})\rho\,dx,$$ where the first one is the momentum of
inertia.

We consider below the solutions to \eref{1.1} -- \eref{1.3} such
that the integrals $m,P, E$ and $G$ converge and call them {\it
solutions with finite momentum of inertia} \,(FMI). It is easy to
verify that for this class of solutions
$M(t)=M=const,\,P(t)=P=const,$
\begin{equation}\label{1.4}
\mathcal E'(t)=-\frac{\Lambda}{\gamma-1}\,\int\limits_{{\mathbb
R}^n} \rho^{1/2}\, p^{3/2}\,dx.
\end{equation}
 The latter equation expresses
the inelastic energy loss per collision.

Further we consider a function $K=p\,\rho^{-\gamma},$ where $\ln K$
can be interpreted as the usual hydrodynamic entropy. System
\eref{1.1} -- \eref{1.3} result
\begin{equation}\label{1.5}
\frac{dK}{dt}=-\Lambda \,K^{\frac{3}{2}}
\,\rho^{\frac{\gamma+1}{2}}\le 0,
$$
\end{equation}
therefore
\begin{equation}\label{1.6}
K(t,x)\le K_+=\sup\limits_{x\in{\mathbb R}^n}K(0,x).
\end{equation}

\section{Main theorem: nonexistence of global smooth solutions}

\begin{theorem}\label{MainT} Let $P\ne 0$ and $M$ be sufficiently small.  Then
there exists no global in time classical FMI solution to the Cauchy
problem for \eref{1.1} -- \eref{1.3}.
\end{theorem}

To prove the theorem we need to get firstly certain estimates of
energy.

\begin{lemma}\label{L1} For the classical FMI solutions to \eref{1.1} -- \eref{1.3} the
following estimates hold:
\begin{equation}\label{2.1}
E_k(t)\ge \frac{P^2}{2M}=const,
\end{equation}
\begin{equation}\label{2.2}
\mathcal E'(t)\le -\Lambda C_1
E_i^{\frac{3\gamma-1}{2(\gamma-1)}}(t),
\end{equation}
where
$C_1=K_+^{-\frac{1}{\gamma-1}}\,(\gamma-1)^{\frac{3\gamma-1}{2(\gamma-1)}}\,
M^{-\frac{\gamma+1}{2(\gamma-1)}}.$
\end{lemma}

{\it Proof.} Inequality \eref{2.1} follows immediately from the
H\"older inequality. To prove \eref{2.2} we firstly use the Jensen
inequality as follows:
$$
\fl \left(\frac{\int\limits_{{\mathbb
R}^n}\,p\,dx}{M}\right)^{\frac{3\gamma-1}{2(\gamma-1)}}=\left(\frac{\int\limits_{{\mathbb
R}^n}\,K\rho\rho^{\gamma-1}\,dx}{M}\right)^{\frac{3\gamma-1}{2(\gamma-1)}}
$$
\begin{equation}\label{2.3}
\fl \le\frac{\int\limits_{{\mathbb R}^n}\,\rho
K^{\frac{3\gamma-1}{2(\gamma-1)}}\,\rho^{\frac{3\gamma-1}{2}}\,dx}{M}
\le\frac{\int\limits_{{\mathbb
R}^n}\,\rho^{1/2}\,p^{3/2}\,K^\frac{1}{\gamma-1} \,dx}{M}\le
\frac{K_+^\frac{1}{\gamma-1}}{M}{\int\limits_{{\mathbb
R}^n}\,\rho^{1/2}\,p^{3/2} \,dx}.
\end{equation}
Together with \eref{1.4} inequality \eref{2.3} gives \eref{2.2}. $
\square$

The following lemmas establish the properties of the momentum of
inertia. Acting as in \cite{Chemin} we get
\begin{lemma}\label{L2}
For classical FMI solutions to \eref{1.1} -- \eref{1.3} the
equalities
\begin{equation}\label{2.4}
G'(t)=F(t),
\end{equation}
\begin{equation}\label{2.5}
F'(t)=2 E_k(t)+n(\gamma-1) E_i(t).
\end{equation}
take place.
\end{lemma}

{\it Proof.} The lemma can be proved by direct calculation using the
general Stokes formula. $\square$

Then we get two-sided estimates of $G(t).$

\begin{lemma}\label{L3}
If $\gamma\le 1+\frac{2}{n}, $ then  for the classical FMI solutions
 to \eref{1.1} -- \eref{1.3} the estimates
 \begin{equation}\label{2.6}
\frac{P^2}{2M}t^2+F(0)t+G(0)\le G(t)\le  \mathcal E(0)
t^2+F(0)t+G(0)
\end{equation}
 hold.
\end{lemma}

{\it Proof.} First of all \eref{2.6}, \eref{2.7} result
\begin{equation}\label{2.7}
G''(t)=2 E_k(t)+n(\gamma-1)E_i (t)=2\mathcal
E(t)-(2-n(\gamma-1))E_i(t).
$$
\end{equation}
Therefore together with \eref{2.1} we have
\begin{equation}\label{2.8}
\frac{P^2}{M}\,\le \, G''(t) \le \,2\mathcal E(0),
\end{equation}
after integration this gives \eref{2.6}. $\square$

Now we get an upper estimate of $E_i(t).$

\begin{lemma} \label{L4} If $\gamma\le 1+\frac{2}{n}$ and $P\ne 0,$ then for
the classical FMI solutions the following estimate is true:
\begin{equation}\label{2.9}
E_i(t)\le \frac{C_2}{G^{n(\gamma-1)/2}},
\end{equation}
where $C_2=\frac{(4G(t)\mathcal
E(t)-F^2(t))G^{(\gamma-1)n/2}(0)}{4}.$
\end{lemma}

{\it Proof.} The method of obtaining the upper estimate of $E_i(t)$
is similar to \cite{Chemin}. Namely, let us consider the function
$Q(t)=4G(t)\mathcal E(t)-F^2(t).$ The H\"older inequality gives
$F^2\le 4 G(t) E_k(t),$ therefore $\mathcal E(t)=E_k(t)+E_i(t)\ge
E_i(t)+\frac{F^2(t)}{4G(t)}$ and
\begin{equation}\label{2.10}
E_i(t)\le \frac {Q(t)}{4 G(t)}.
\end{equation}
We notice also that $Q(t)>0$ provided the pressure does not equal to
zero identically. Then taking into account \eref{2.4}, \eref{2.5}
and \eref{2.7} we have
\begin{eqnarray}\label{2.11}
 Q'(t)=4 G'(t) \mathcal E(t)-2 G'(t) G''(t)+ 4G(t)\mathcal
E'(t)\nonumber\\= 2(2-n(\gamma-1))\,G'(t)E_i(t)+4G(t)\mathcal E'(t).
\end{eqnarray}
Further, one can see from \eref{2.8} that  $G'(t)>0$ beginning from
a positive $t_0$ for all initial data. Thus, for $\gamma\le
1+\frac{2}{n}$ we have from \eref{2.10}, \eref{2.11}
\begin{equation}\label{2.12}
\frac{Q'(t)}{Q(t)}\le \frac{2-n(\gamma-1)}{2}\,\frac{G'(t)}{G(t)}.
\end{equation}
Then \eref{2.10} and \eref{2.12} give
$$E_i(t)\le \frac{C_2}{G^{(\gamma-1)n/2}(t)}, \quad C_2=\frac{Q(0)G^{(\gamma-1)n/2}(0)}{4}.$$
The proof is over. $\square$

\begin{remark}
Inequalities \eref{2.1} and \eref{2.2} result
$$
\mathcal E'(t)\le -\Lambda \,C_1\,\left(\mathcal
E(t)-\frac{P^2}{2M}\right)^{\frac{3\gamma-1}{2(\gamma-1)}}.
$$
Integrating this inequality we the following upper estimate of
$E_i(t):$
$$E_i(t)\le (c_1\,t+c_2)^{-\frac{2(\gamma-1)}{\gamma+1}},$$
where $c_1=\frac{\Lambda C_1(\gamma+1)}{\gamma-1}$ and
$c_2=(E(0)-P^2/(2M))^{\frac{2(\gamma-1)}{\gamma+1}} $ are positive
constants. However, for $\gamma \le 1+\frac{2}{n}$ this estimate is
less exact than \eref{2.9}  and is not enough for our proof.
\end{remark}

The next step is a lower estimate of $E_i(t).$

\begin{lemma} \label{L5} Let $P\ne 0.$ Then for the classical FMI solutions the estimate
\begin{equation}\label{2.13}
E_i(t)\ge \frac{C_3}{G^{(\gamma-1)n/2}(t)},
\end{equation}
holds with a positive constant
$$C_3=\frac{1}{2^{(\gamma-1)n/2}\,K_+^{(\gamma-1)(1+n/2)}\,
(\gamma-1)}\,{C_{\gamma,n}^{-\frac{\gamma(n+2)-2}{2}}}\,
(C_6(\gamma)) ^{\frac{(\gamma-1)\gamma(n+2)-2}{4\gamma-1}},$$ the
value of $C_6$ is written in \eref{2.19}.
\end{lemma}

{\it Proof.}
 The proof is based on the inequality
$$
\|f\|_{L^1({\mathbb
R}^n;\,dx)}\,\le\,C_{\gamma,n}\,\|f\|^{\frac{2\gamma}{(n+2)\gamma-n}}_{L^\gamma({\mathbb
R}^n;\,dx)}\,\|f\|^{\frac{n(\gamma-1)}{(n+2)\gamma-n}}_{L^1({\mathbb
R}^n;\,|x|^2\,dx)},
$$
$$C_{\gamma,
n}=\left(\frac{2\gamma}{n(\gamma-1)}\right)^
{\frac{n(\gamma-1)}{(n+2)\gamma-n}} +
\left(\frac{2\gamma}{n(\gamma-1)}\right)^
{\frac{-2\gamma}{(n+2)\gamma-n}},$$ established in \cite{Chemin}.

Namely, we have for $f=K\rho$
$$
E_i(t)=\frac{1}{\gamma-1}\,\int\limits_{{\mathbb
R}^n}\,\frac{p}{\rho^\gamma} \,\frac{(K\rho)^\gamma}{K^\gamma}\,dx
\ge \frac{1}{K_+^{\gamma-1}\,(\gamma-1)}\,\int\limits_{{\mathbb
R}^n}\, (K\rho)^\gamma\,dx\,$$
$$
\ge\,\frac{1}{K_+^{\gamma-1}\,(\gamma-1)}\,\frac{\left(C^{-1}_{\gamma,n}\,\int\limits_{{\mathbb
R}^n}\,
K\rho\,dx\right)^{\frac{\gamma(n+2)-2}{2}}}{\left(\int\limits_{{\mathbb
R}^n}\, K\rho |x|^2\,dx\right)^{\frac{n(\gamma-1)}{2}}}$$
\begin{equation}\label{2.14}
\ge
\,\frac{1}{2^{(\gamma-1)n/2}\,K_+^{(\gamma-1)(1+n/2)}\,(\gamma-1)}\,\frac{\left(C^{-1}_{\gamma,n}\,S(t)
\right)^{\frac{\gamma(n+2)-2}{2}}}{(G(t))^{\frac{n(\gamma-1)}{2}}},
\end{equation}
where we denoted $S(t)=\int\limits_{{\mathbb R}^n}\, K\rho\,dx.$

Further, from \eref{1.1} and \eref{1.5} we have
$$
\partial_t(K\rho)+{\rm div}_x (K\rho u)=-\Lambda K^{\frac{3}{2}}
\rho^{\frac{\gamma+3}{2}},
$$
therefore
\begin{equation}\label{2.15}
S'(t)=-\Lambda\,\int\limits_{{\mathbb R}^n}\,
K^{\frac{3}{2}}\,\rho^{\frac{\gamma+3}{2}}\,dx.
\end{equation}
 From
the Jensen inequality we get
$$
\int\limits_{{\mathbb R}^n}\,
K^{\frac{3}{2}}\,\rho^{\frac{\gamma+3}{2}}\,dx\le
\frac{K_+^{\frac{3}{2}}\,M}{S(t)}\int\limits_{{\mathbb R}^n}\,
K\,\rho\,\rho^{\frac{\gamma+1}{2}}\,dx\le K_+^{\frac{3}{2}} M
\,\left(\frac{\int\limits_{{\mathbb R}^n}\,
K\,\rho^{\gamma}\,dx}{S(t)} \right)^{\frac{1+\gamma}{2(\gamma-1)}}$$
\begin{equation}\label{2.16}
= C_4 \,\left(\frac{E_i(t)}{S(t)}
\right)^{\frac{1+\gamma}{2(\gamma-1)}},
\end{equation}
where $C_4=K_+^{\frac{3}{2}}
M(\gamma-1)^{\frac{1+\gamma}{2(\gamma-1)}}.$ Further, from
\eref{2.9}, \eref{2.15}, \eref{2.16} we obtain
\begin{equation}\label{2.17}
S^{\frac{1+\gamma}{2(\gamma-1)}}(t)\,S'(t)\ge -\Lambda
C_4\,C_2^{\frac{\gamma+1}{2(\gamma-1)}}\,(G(t))^{-\frac{n(\gamma+1)}{4}}
.
\end{equation}
Now we take into account the lower estimate in \eref{2.6} together
with the fact that beginning from a certain $t_0$ the value of
$F(t)$ becomes positive if $P\ne 0$ (see \eref{2.1}, \eref{2.5}) and
integrate \eref{2.17}. Thus for $t\ge t_0$ we get in the case
$\gamma<3$
\begin{equation}\label{2.18}
S(t)\ge\left(S(t_0)^{\frac{\gamma-3}{2(\gamma-1)}}+
\Lambda\frac{4(\gamma-1)}{(3-\gamma)(n(\gamma-1)-2)}\,C_5^{1-\frac{n(\gamma+1)}{2}}\right)^{\frac{2(\gamma-1)}{\gamma-3}},
\end{equation}
where $C_5=M\,(F(t_0))^{-\frac{n(\gamma+1)}{2}}\,P^{-2},$ and for
$\gamma=3$ \,$(n=1)$
\begin{equation}\label{2.19}
S(t)\ge S(t_0)\,\exp \left(\frac{2\Lambda}{C_5} \right).
\end{equation}
We denote the constant in the right-hand side of \eref{2.18},
\eref{2.19} by $C_6(\gamma).$

It is easy to see that for sufficiently large  $t$ the value of
$S(t)$ is separated from zero.

Thus, from \eref{2.14}, \eref{2.18}, \eref{2.19} we obtain
\eref{2.13}. $\square$

\subsubsection*{Proof of Theorem \ref{MainT}}

Taking into account \eref{2.2}, \eref{2.13} and the right-hand side
of \eref{2.6} we get
\begin{eqnarray}\label{2.20}
 \mathcal E'(t)\le -\Lambda C_1
\,C_3^{\frac{3\gamma-1}{2(\gamma-1)}}\,(G(t))^{-\frac{n(3\gamma-1)}{4}}\\\nonumber
\le -\Lambda C_1 \,C_3^{\frac{3\gamma-1}{2(\gamma-1)}}\,(\mathcal E
(0) t^2 +F(0)t+G(0))^{-\frac{n(3\gamma-1)}{4}}.
\end{eqnarray}
As follows from \eref{2.1}, \eref{2.4} beginning from $t_0>0$ the
value of $F(t)$ becomes positive.  Integrating \eref{2.20} we obtain
for $t>t_0$
$$
\mathcal E(t)\le \mathcal E(t_0)-
\Lambda\,C_1\,C_3^{\frac{3\gamma-1}{2(\gamma-1)}}\,
\frac{2^{n(3\gamma-1)/2}(F(t_0))^{-n(3\gamma-1)/2+1}(\mathcal
E(0))^{n(\gamma-1)/4-1}}{n(3\gamma-1)-2}+\lambda(t),
$$
where $\lambda(t)=\Or(t^{-n(3\gamma-1)/2+1}),\,t\to\infty.$ Since
$C_3$ tends to a positive constant as $M\to 0$ and $C_1$ contains
$M$ in a negative degree (see the statements of Lemmas \ref{L1} and
\ref{L5}, then choosing $M$ sufficiently small, we can always get a
contradiction with inequality \eref{2.1}. This prove the theorem.
$\square$

\begin{remark}
Main idea of this paper can be found in \cite{Rozanova}, where it
was proved nonexistence of global smooth solutions to the
compressible Navier-Stokes equations. For the freely cooling gas, it
would be possible to add viscous terms (with constant viscosity
coefficients) to apply the same technique and obtain the analogous
nonexistence result.
\end{remark}

\section{One-dimensional case}

As it was noticed system \eref{1.1} -- \eref{1.3} has no constant
solution except of the trivial one $p \equiv 0.$ However, it is
possible to construct the nontrivial steady state  solution
$\bar\rho(0,x), \bar v(0,x), \bar p(0,x)$ for the regions
$|x|>R(t)>0.$  For $|x|\le R(0)$ we chose the functions
$\,\rho(0,x), \,v(0,x), \,p(0,x)\,$ arbitrarily to get initial data,
smooth on the whole real axis.  We are going to show that such
solution necessarily loses its initial smoothness.

\subsection{Automodel solution}

Let us find a solution that depends on the automodel variable
$\xi=x-at,\,a=const.$ The continuity equation \eref{1.1} gives the
connection between the velocity and density as follows:
\begin{equation}\label{3.1}
u(\xi)=\frac{c_1}{\rho(\xi)}+a, \, c_1=const\ne 0.
\end{equation}
Equation \eref{1.2} and \eref{3.1} results
\begin{equation}\label{3.2}
c_1 u(\xi)+p(\xi)=c_2,\,c_2=const.
\end{equation}
From \eref{3.1} and \eref{3.2} we have $c_1^2+\rho p=\rho (c_2-a
c_1)$, therefore $c_2-a c_1>0$.

 Further, we substitute the functions $u$ and $p,$ found from
(3.1) and (3.2) and expressed through $\rho,$  in the equation
$$
\partial_t p + u\partial_x p + \gamma p \partial_x u = -\Lambda
\rho^2 p^{3/2},
$$
which is a corollary of \eref{1.1}, \eref{1.3} and the state
equation $p=\rho T$. Thus we get an ordinary differential equation
$$
\rho'(\xi)=\,-\,\frac{\Lambda}{c_1}\,\frac{\rho^2(\xi)\,((c_2-ac_1)\rho(\xi)-c_1^2
)^{3/2}}{c_1^2(\gamma+1)-\gamma(c_2-ac_1)\rho(\xi)},
$$
or
\begin{equation}\label{3.3}
z'(\xi)=-\frac{\Lambda}{c_1(c_2-ac_1)}\,\frac{z^{3/2}(\xi)(z(\xi)+c_1^2)^2}{c_1^2-\gamma
z(\xi)},
\end{equation}
where $z=(c_2-ac_1)\rho-c_1^2. $ The case $c_2-ac_1=0 $ that seems
simpler corresponds to a negative pressure and we do not consider
it. Equation \eref{3.3} has a solution
\begin{equation}\label{3.4}
\fl
\xi=f(z):=\frac{c_2-ac_1}{c_1\,\Lambda}\left(\frac{(\gamma+3)\,\arctan\frac{\sqrt{z}}{c_1}}{
c_1}\,+\, \,\frac{(\gamma+1)\sqrt{z}}{z+c^2_1}\,+\,
\frac{2}{\sqrt{z}}\right)\,+\, c_3,
\end{equation}
$\,c_3=const$. One can see from \eref{3.3} that if $z(\xi_0)\le
z_*=\frac{c_1^2}{\gamma}$, then there exists a branch of solution
\eref{3.4} defined on the semi-axis  $(\xi_0,+\infty)$  for positive
$c_1\,$ where the function $z(\xi)$ decrease monotonically from the
value $z(\xi_0)$  to zero. For negative $c_1\,$ this branch is
defined  on the semi-axis $(-\infty, \xi_0),$ where the function
$z(\xi)$ increase monotonically from zero to $z(\xi_0)$.


\subsection{"Finally steady state" and its smooth compact perturbation}

Let us consider a stationary solution for $a=0,$ $\xi=x$,  choose a
point $x_0\ge 0$, a constant $c_1=u(x_+)\rho(x_+)>0$,
and construct on the semi-axis $x>x_+ $ a solution $z_+(x)$.
Analogously for the semi-axis $x\le x_- $ we choose
$c_1=u(x_-)\rho(x_-)<0$ and construct a solution $z_-(x).$ Thus,
outside of the segment $[x_-, x_+]$ we define a solution
$$\bar z(x)=\left\{\begin{array}{ll} z_-(x),&\quad x\in (-\infty,\, x_-), \\ z_+(x), &\quad x\in
(x_+, +\infty).
\end{array}\right.$$

For the sake of simplicity we set $x_-=-x_+$, $c_1=-k$ for $z_-(x)$
and $c_1=k $ for $z_+(x)$, where $k=const>0$.

In their turn, the density, pressure and velocity can be found as
\begin{equation}\label{3.5}
$$\bar\rho=\frac{\bar z+k^2}{c_2},\quad \bar p=c_2-\frac{k^2}{\bar \rho},\quad \bar u=\frac{k\,{\rm sign} x
\,}{\bar \rho}. \end{equation}

It is very attractive to choose $x_+=x_-=0$ and to construct a
piecewise continuous solution like a solution of a "Riemmann
problem" with non-constant left and right states. Nevertheless, it
can be readily shown that the Hugoniot conditions do not hold on the
jump. Indeed, the components $\bar\rho$ and $\bar p$ are continuous
at the point $x=0,$ having a jump in derivative, however, the
velocity itself has a jump $[\bar u]=\frac{2k}{\bar \rho(0)}\ge
\frac{2\gamma c_2}{k(\gamma+1)}$ (see Figs.1 -- 3). The value of
$[\bar u]$ tends to zero as $k\to\infty,$ however the Hugoniot
condition $[\bar \rho \bar u]=0$ does not implement in the origin
$x=0.$ Of course, we can choose $c_1$ such that the density and
pressure have jumps in the origin and consider $a\ne 0$,
nevertheless the careful analysis shows that the Hugoniot conditions
 do not hold anyway.

\begin{figure}[h]
\centerline{\includegraphics[width=0.7\columnwidth]{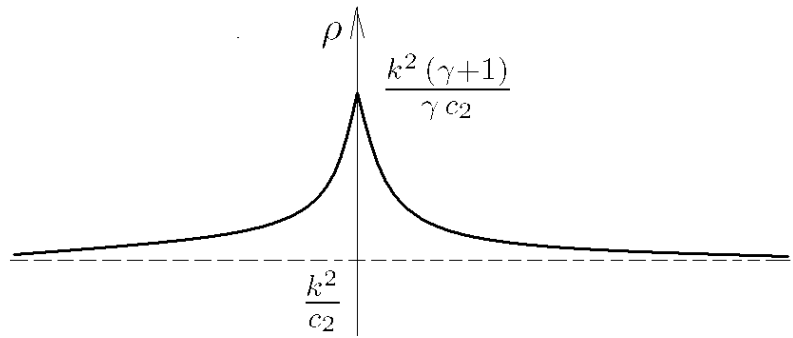}}
\caption{The density}
\end{figure}%
\begin{figure}[h]
\centerline{\includegraphics[width=0.7\columnwidth]{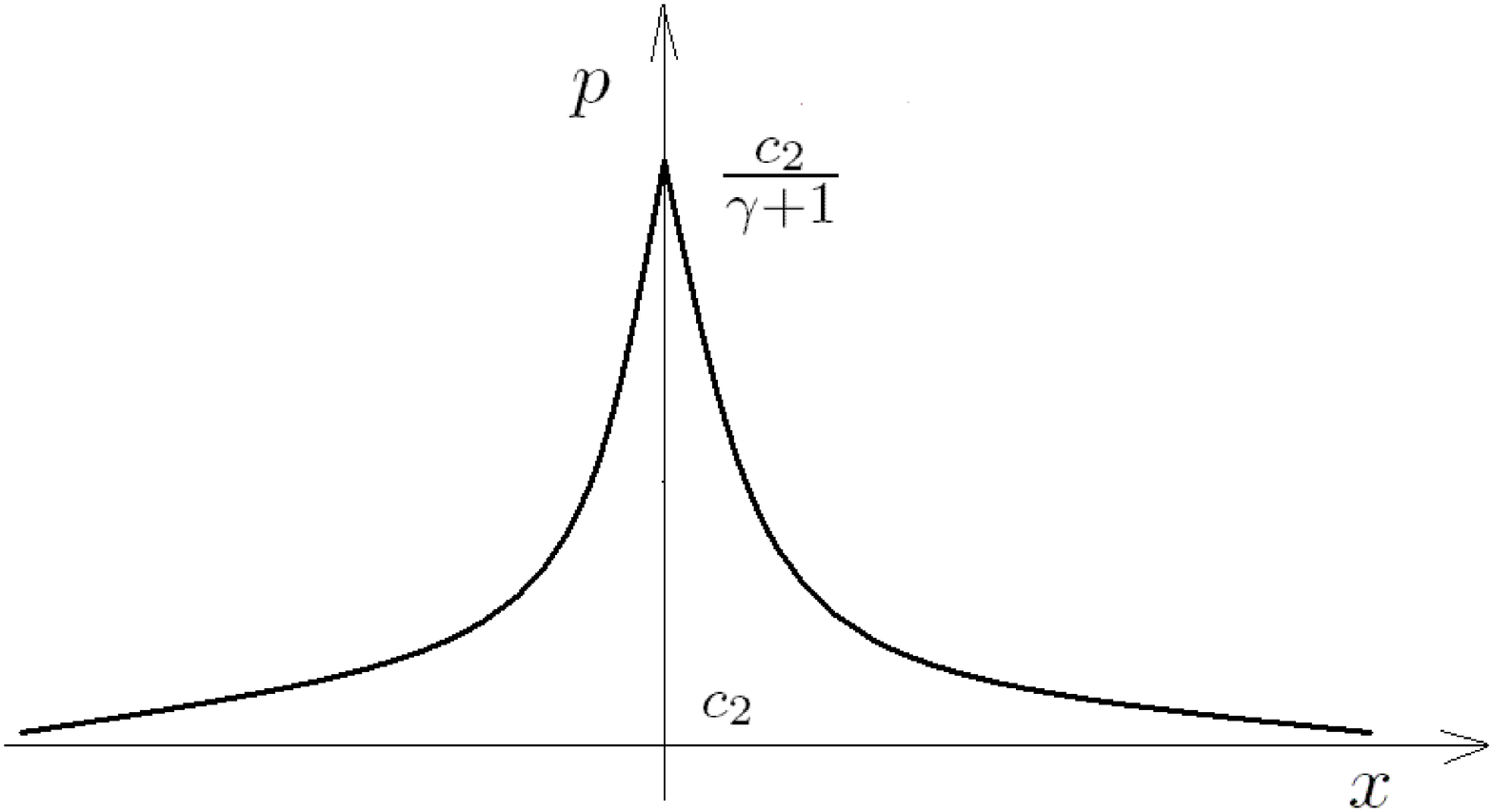}}
\caption{The pressure}
\end{figure}%
\begin{figure}[h]
\centerline{\includegraphics[width=0.7\columnwidth]{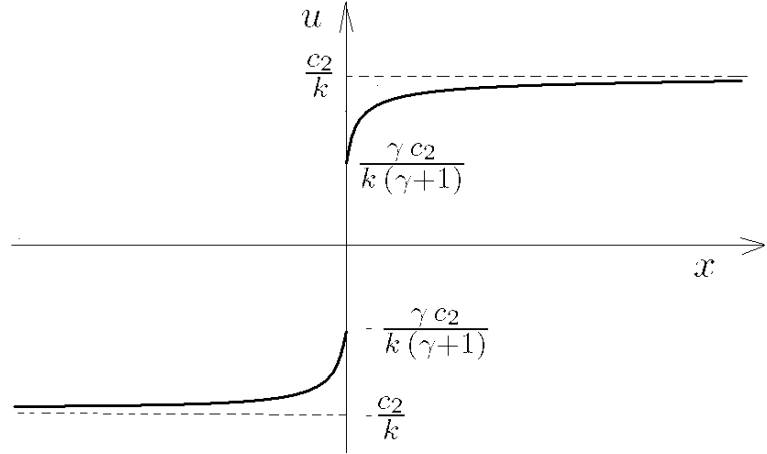}}
\caption{The velocity}
\end{figure}%


Therefore we choose $R(0)>x_+,$ and define smooth initial data such
that for $|x|>R(0)$ they coincide with $\bar \rho(x),\, \bar u(x),
\,\bar p(x)$ and on the segment $|x|\le R(0)$ they are arbitrary
smooth functions $\,\rho(0,x)>0, \, p(0,x)>0,\, u(0,x). $ We will
call this type of initial data {\it the compact perturbation of
nontrivial finally steady state.} Let us note that for $k=0$ we get
the trivial zero-state solution.

\subsection{Breakdown of the compact perturbation of the nontrivial "finally steady
state"}\label{SecFSS}

Let us denote the perturbed region by $B(t):=\{x\big| \,|x|\le R(t)
\},$ and consider the  analogs of functionals used in the previous
section:
$$\tilde
G(t)=\frac{1}{2}\,\int\limits_{B(t)}\, \rho (t,x) |x|^2\,dx,\qquad
\tilde F(t)=\int\limits_{B(t)}\, \rho (t,x) (u(t,x),x)\,dx,
$$
$$\tilde
E_k(t)=\frac{1}{2}\,\int\limits_{B(t)}\, \rho (t,x) |u(t,x)|^2\,dx,
\qquad \tilde M(t)=\int\limits_{B(t)}\, \rho (t,x) \,dx.
$$
$$\tilde
P(t)=\,\int\limits_{B(t)}\, \rho (t,x)\,u(t,x) \,dx,\quad \tilde
S(t)=\,\int\limits_{B(t)}\,K(t,x) \rho (t,x) \,dx.
$$

The following theorem holds:

\begin{theorem} Let $\tilde P^2
(0)>\frac{8\tilde M c_2}{k},\,n=1.$ Then there exists no globally in
$t$ smooth perturbation of the nontrivial steady state for system
\eref{1.1}-- \eref{1.3}.
\end{theorem}

{\it Proof.} First of all we note that  system \eref{1.1}--
\eref{1.3} is hyperbolic and therefore the speed of boundary of the
perturbations is equals to $ |\bar u|+ V_s,$ where
$V_s=\sqrt{\bar
p_{\bar \rho}},$ the sound speed. As follows from \eref{3.5},
\begin{equation}\label{3.6}
V_s=\frac{k}{\bar \rho}\le \frac{ c_2}{k}:=\sigma,\end{equation}
where we use the estimate $\bar \rho \ge \frac{k^2}{c_2}.$

Then we have
\begin{equation}\label{3.7}
\tilde F'(t)=2 \tilde E_k + \int\limits_{B(t)}\, (p (t,x)-\bar
p(t,R(t))) \,dx,
\end{equation}
$$\tilde M'(t)=0,\quad \tilde P'(t)=0.$$
Further,  the H\"older inequality implies
\begin{equation}\label{3.8}
\tilde E_k(t)\ge \frac{\tilde F^2(t)}{4 \tilde G(t)},
\end{equation}
and one can estimate
\begin{equation}\label{3.9}
\tilde G(t)\le \frac{1}{2}\,R^2(t)\tilde M(t) \le
\frac{1}{2}\,R^2(t)\tilde M(0).
\end{equation}
Further, the Jensen inequality yields
\begin{equation}\label{3.10}
\int\limits_{B(t)}\, p (t,x) \,dx\ge (2
R(t))^{1-\gamma}\,K_+^{\gamma-1} \,(\tilde S(t))^\gamma.
\end{equation}
As in the proof of Lemma \ref{L5} we can show that for sufficiently
large $t$ the function $S(t)\ge S_-=const >0.$ Then we notice that
\eref{3.3} implies
\begin{equation}\label{3.11}
\bar z(x)\sim \frac{4k^2 c^2_2}{\Lambda^2}\,x^{-2}, \quad
x\to\infty.
\end{equation}
Therefore taking into account \eref{3.5} we have
\begin{equation}\label{3.12}
\int\limits_{B(t)}\, \bar p(t,R(t))) \,dx\,\le
\,\frac{c_2}{k^2}\,\int\limits_{B(t)}\, \bar z(R(t))
\,dx\,\sim\,\frac{4c_2^2}{\Lambda^2}\,R^{-1}(t), \quad t\to\infty.
\end{equation}
Thus, from \eref{3.6} -- \eref{3.12} and the estimate $$\tilde
E_k(t)\ge \frac{\tilde P(t)}{2\tilde M(t)}=\frac{\tilde
P(0)}{2\tilde M (0)}$$ we have
\begin{eqnarray}\label{3.13}\fl
\tilde F'(t)\,\ge\,\tilde E_k(t)\,+\,\frac{\tilde F^2(t)}{4\, \tilde
M (0)\,(R(0)+\sigma
t)^2}\,+\,\frac{K_+^{\gamma-1}\,S_-^\gamma}{(R(0)+\sigma
t)^{\gamma-1}}\,-\,\frac{4\,c^2_2}{\Lambda^2\,(R(0)+\sigma
t)}\,\\\nonumber \ge \,\frac{\tilde P(0)}{2\tilde M(0)}\, +
\,\frac{\tilde F^2(t)}{4\, \tilde M (0)\,(R(0)+\sigma
t)^2}\,+\,\lambda(t),
\end{eqnarray}
where $\lambda(t)\to 0,\,t\to\infty.$

As we can see from \eref{3.13}, beginning from a certain $t_0>0$ the
function $\tilde F(t)>4\,\tilde M(0)\,R(t_0)\,\sigma,$ moreover,
\begin{equation}\label{3.13a}
\tilde F'(t)\,\ge\,\frac{\tilde F^2(t)}{4\, \tilde
M(0)\,(R(0)+\sigma t)^2},\quad t>t_0.
\end{equation}
 Integrating \eref{3.13a}
 from $t=t_0$ we get
\begin{equation*}\nonumber
\tilde F(t)\,\ge \,\frac {4\,\tilde M(0) \tilde F(t_0)
R(t_0)(R(t_0)+\sigma t)}{4\tilde M(0) R^2(t_0)+(4\tilde M(0)
R(t_0)\sigma -\tilde F(t_0))t}.
\end{equation*}
Thus, $\tilde F(t)$ blow-ups at a finite time. This contradicts to
the inequality $\tilde F^2(t)\le 4 \tilde G(t) \tilde E(0).$ The
theorem is proved. $\square$

\begin{remark}
The idea  of the method is due to \cite {Sideris}, where it was
proved that the compact smooth perturbation of a constant state of
gas dynamics equation can not be globally smooth in time.
\end{remark}

\section{Exact solution with singularity}

 Naturally there arises a question on a type of
predicted singularity. In particular, in the remarkable papers
\cite{Meerson1}, \cite{Meerson2} for the one-dimensional case the
authors employ Lagrangian coordinates and derive a broad family of
exact non-stationary non-self-similar solutions. These solutions
exhibit a singularity, where the density blowups in a finite time
when starting from smooth initial conditions. Moreover, the velocity
gradient also blowups while the velocity itself and develop a cusp
discontinuity (rather then a shock) at the point of singularity.
This approach is partially extended  to the 2D case in
\cite{Fouxon}.

Here for any spatial dimensions we construct a simple family of
solutions to the system \eref{1.1} -- \eref{1.3} having a
singularity in the density whereas other components are continuous.
Indeed, if we substitute in \eref{1.1} -- \eref{1.3}
\begin{equation}
\label{3.14} u(t,x)=\alpha(t)\,x,\quad
\rho(t,x)=\beta(t)\,|x|^q,\quad p(t,x)=s(t)\,|x|^l, \end{equation}
where $x$ is a radius-vector of point, we obtain
\begin{equation}
\label{3.15}  q=-1,\qquad l=1,\qquad\beta(t)=\beta_0=const\ge 0,
\end{equation} and $\alpha(t),\,s(t)\ge 0$ satisfy the following
system of nonlinear ODE:
\begin{equation}\label{3.16}
\alpha'(t)+\alpha^2(t)+\frac{s(t)}{\beta_0}=0,
\end{equation}
\begin{equation}\label{3.18}
s'(t)+(\gamma+1)\,n\,
\,s(t)\,\alpha(t)=-\Lambda\,\beta^{1/2}(t)\,s^{3/2}(t).
\end{equation}

This system has a unique equilibrium $(\alpha(t)=0,\, s(t)=0)$, it
is unstable. One of its solutions is very simple: $s(t)=0\, (p\equiv
0),\, \alpha(t) =(t+\alpha^{-1}(0))^{-1}$. An analysis of the phase
portrait shows that if $s(0)>0$, then $\alpha(t)\to -\infty,\,
s(t)\to +\infty$ for all $\alpha(0).$

Let us prove this fact in a different way. We consider a symmetric
material volume $B(t) $ containing the origin $x=0$ and use the
denotation of Sec.\ref{SecFSS}. We can see that in spite of the
singularity in the component of density all integrals below exist.
Thus, due to the structure of solution \eref{3.14}, \eref{3.15} we
have
\begin{equation}\label{3.19}\fl
\tilde G'(t)=\tilde F(t)=2\alpha(t)\,\tilde G(t),
\end{equation}
\begin{equation}\label{3.20}
\fl \tilde F'(t)=2\,\tilde E_k(t)+\int\limits_{B(t)}\, (p (t,x)-\bar
p(t,R(t))) \,dx\le 2\,\tilde E_k(t)=2\,\alpha^2 (t)\,\tilde G(t).
\end{equation}
As follows from \eref{3.19}, \eref{3.20}, the velocity gradient
obeys the inequality
$$
\alpha'(t)\le -\alpha^2(t),
$$
therefore $\alpha(t)\le ({t+\alpha^{-1}(t_0)})^{-1}, $ and in the
case $\alpha(t_0)<0$ we can see that  $\alpha(t)\to -\infty $ as
$t\to
 -\alpha^{-1}(t_0).$ Since \eref{3.16} -- \eref{3.18}
 result that $\alpha'(t)< -{s(t)}/{\beta_0}\le - \varepsilon <0$
 (the latter inequality follows from the uniqueness theorem),
 $\alpha(t)$ become negative in a finite time. The proof of over.

We see that in the presence of a stationary singularity in the
component of density, a balance between velocity and pressure arise.
Generically the velocity and pressure blow up in a finite time.
Thus, the components of this "black hole" solution solution can
collapse in different moments of time. The singularity of density in
the origin is integrable for the dimension $n\ge 2$.

Let us remark that for usual gas dynamics the solutions of such kind
with a ``linear profile'' of  velocity  are well investigated (e.g.
\cite{Sedov}, Chapter IV, Sec.15).

The author is  indebted  to  B.Meerson for attracting attention to
the problem and thanks N.Leontiev and V.Shelkovich for a helpful
discussion.


\section*{References}

\begin{thebibliography}{99}

\bibitem{Brilliantov} Brilliantov N V and  P\"oschel T 2004. {\em Kinetic theory of granular
gases}  (Oxford: Oxford University Press)

\bibitem{Chemin}Chemin J-Y \, 1990 \, Dynamique  des  gaz  \`a masse
totale finie,\,
  {\em  Asymptotic Analysis} {\bf {3}} 215-20

\bibitem{Meerson1}  Fouxion I,  Meerson B,   Assaf M and Livne E\, 2007\,
 Formation and evolution of density singularities in
hydrodynamics of inelastic gases  {\em Phys.Rev. E} {\bf 75}, 050301
(R)

\bibitem{Meerson2} Fouxion I,  Meerson B,   Assaf M and Livne E\, 2007\, Formation
and evolution of density singularities in ideal hydrodynamics of
freely cooling inelastic gases: A family of exact solutions {\em
Physics of Fluids}  {\bf 19}, 093303

\bibitem{Fouxon}  Fouxon I \,200, Finite-time collapse and localized states in the dynamics of
dissipative gases {\em Phys. Rev. E}{\bf 80} 010301(R)

\bibitem{Kato} Kato T \, 1975 \, The Cauchy problem for quasilinear
symmetric hyperbolic systems {\em Arch.Ration.Math.Anal.}{\bf 58}
181-205

\bibitem{Ludvig}  Luding S 2009  Towards dense, realistic granular media in 2D
{\em Nonlinearity} {\bf 22} R101

\bibitem{Rozanova} Rozanova O \, 2008\, Blow up of smooth highly decreasing at
infinity solutions to the compressible Navier-Stokes equations {\em
Journal of Differential Equations} {\bf{245}} 1762-74

\bibitem{Sedov}  Sedov L I \,1982\, Similarity and dimensional methods in mechanics
(Moscow: Mir)

\bibitem{Sideris}  Sideris T C \,1985 \,Formation of singularities in three-dimensional
compressible fluids {\em Comm.Math.Phys.} {\bf 101} 475-85





\end{thebibliography}
\end{document}